\documentclass[leqno,11pt,a4paper]{article}
\usepackage[left=3.2cm, top=2.2cm,bottom=2.5cm]{geometry}
\usepackage{amsmath,amssymb,amsthm,mathrsfs,calc,graphicx}
\usepackage[british]{babel}

\numberwithin{equation}{section}

\newtheorem{theorem}{Theorem}[section]
\newtheorem{proposition}{Proposition}[section]
\newtheorem{lemma}{Lemma}[section]

{\theoremstyle{definition}

}
{\theoremstyle{theorem}
\newtheorem {remark}{Remark}[section]

}

\newcommand{\dist}{\operatorname{dist}}
\newcommand{\diam}{\operatorname{diam}}
\newcommand{\RR}{\mathbb{R}}
\newcommand{\EE}{\mathbb{E}}
\newcommand{\NN}{\mathbb{N}}
\newcommand{\PP}{\mathbb{P}}
\newcommand{\VV}{\mathbb{V}}
\newcommand{\QQ}{\mathbb{Q}}
\newcommand{\LG}{\mathbb{G}}
\newcommand{\AG}{\mathbb{A}}

\newcommand{\cN}{\mathcal{N}}

\newcommand{\dint}{{\rm d}}

\begin{document}

\title{\bfseries Distances between Poisson $k$-flats}

\author{Matthias Schulte\footnotemark[1]\, and Christoph Th\"ale\footnotemark[2]}

\date{}
\renewcommand{\thefootnote}{\fnsymbol{footnote}}
\footnotetext[1]{Universit\"at Osnabr\"uck, Institut f\"ur Mathematik, Albrechtstr. 28a, D-49076 Osnabr\"uck, Germany. Email: matthias.schulte@uni-osnabrueck.de}

\footnotetext[2]{Universit\"at Osnabr\"uck, Institut f\"ur Mathematik, Albrechtstr. 28a, D-49076 Osnabr\"uck, Germany. Email: christoph.thaele@uni-osnabrueck.de}

\maketitle

\begin{abstract}
The distances between flats of a Poisson $k$-flat process in the $d$-dimensional Euclidean space with $k<d/2$ are discussed. Continuing an approach originally due to Rolf Schneider, the number of pairs of flats having distance less than a given threshold and midpoint in a fixed compact and convex set is considered. For a family of increasing convex subsets, the asymptotic variance is computed and a central limit theorem with an explicit rate of convergence is proven. Moreover, the asymptotic distribution of the $m$-th smallest distance between two flats is investigated and it is shown that the ordered distances form asymptotically after suitable rescaling an inhomogeneous Poisson point process on the positive real axis. A similar result with a homogeneous limiting process is derived for distances around a fixed, strictly positive value. Our proofs rely on recent findings based on the Wiener-It\^o chaos decomposition and the Malliavin-Stein method.
\end{abstract}

\section{Introduction}\label{intro}

Point processes of $k$-dimensional flats in $\RR^d$, especially Poisson point processes, are one of the most classical topics considered in stochastic geometry; cf. \cite{Mecke91,MeckeThomas,Weil87} for early works, \cite{BL,HugLastWeil,Spodarev01,Spodarev03} for more recent papers and the book \cite{SW} for an exhaustive reference. One of the problems considered in the theory of (Poisson) $k$-flat processes, the so-called \textbf{proximity problem}, is to describe the closeness or denseness of the arrangement of the flats in the case $k<d/2$, where the flats do not intersect each other (at least under suitable additional assumptions on their distribution). The notion of proximity generalizes the well-known second-order intersection density for $k$-flat processes in $\RR^d$ with $k\geq d/2$ to the case $k<d/2$ and was originally introduced by Schneider in \cite{Schneider99}. There, only mean values of the proximity functional were considered, but no higher-order moments, limit theory or extreme values.

In this paper, we focus our attention to the Poisson case, for which we compute the asymptotic variance of the classical proximity as considered in \cite{Schneider99} and establish a Berry-Esseen-type central limit theorem. Moreover, we will not only deal with a cumulative proximity functional, but also investigate the order statistics induced by all distances between pairs of distinct flats, in particular the minimal distance, and the behaviour of the distances around a given positive value. This alternative approach to the proximity problem gives new insight into the geometry of Poisson $k$-flat processes.

The proofs of our limit theorems make use of a general central limit theorem from \cite{Schulte2012} and a result about point process convergence and extreme value theory in \cite{SchulteThaele12}. They are based on Berry-Esseen type inequalities in \cite{Peccati12,PTSU10} that were derived by combining the Malliavin calculus of variations on the Poisson space with Stein's method. The backbone of these methods is the fact that each square integrable Poisson functional can be represented as orthogonal sum of multiple Wiener-It\^o integrals; see \cite{LP} and the references therein. It has recently turned out that this so-called Wiener-It\^o chaos decomposition and related limit theorems can successfully be applied to problems in stochastic geometry. For example, in \cite{RS12} a general set-up was investigated as well as central limit theorems for Poisson hyperplanes, \cite{LPST} deals with moment formulas and very general geometric functionals of Poisson $k$-flat processes, \cite{LRPeccati1,LRPeccati2} consider fine Gaussian fluctuations on the Poisson space and geometric random graphs. In all these works a crucial r\^ole is played by a special class of Poisson functionals, the so-called Poisson U-statistics.

The text is structured as follows: In Section \ref{sec:results}, we introduce the proximity of a Poisson $k$-flat process and present our main results, Theorem \ref{thm:Expectation} -- \ref{thm:Maximum}. Their proofs rely on the Wiener-It\^o chaos decomposition of Poisson functionals, whose background is briefly introduced in Section \ref{sec:chaos}. The remaining three sections are devoted to the detailed proofs of our theorems.

\section{Statement of the main results}\label{sec:results}

\subsection{Framework}
A Poisson process of $k$-dimensional flats in $\RR^d$ is a Poisson point process on the space $\AG_k^d$ of $k$-dimensional affine subspaces of $\RR^d$, where $k\in\{1,2,\ldots,d-1\}$ and $d\geq 1$. We let $\eta_t$ be such a \textbf{Poisson process of $k$-flats} having its \textbf{intensity measure} $\Theta_t$ given by
\begin{equation}\label{eq:intensitymeasure}
\int_{\AG_k^d}f(E)\,\Theta_t(\dint E)=t\int_{\LG_k^d}\int_{L^\perp}f(L+x)\,\ell_{E^\perp}(\dint x)\,\QQ(\dint E).
\end{equation}
Here, $f:\AG_k^d\rightarrow\RR$ is a non-negative measurable function, $t>0$, $\LG_k^d$ is the Grassmannian of $k$-dimensional linear subspaces of $\RR^d$, $\ell_{E^\perp}$ is the Lebesgue measure on $E^\perp$ and $\QQ$ is a probability measure on $\LG_k^d$. The Poisson $k$-flat process $\eta_t$ is {\bf stationary}, i.e., its distribution is invariant under all translations. In case that $\QQ$ is the invariant probability measure (Haar measure) $\nu_k$ on $\LG^d_k$, the distribution of $\eta_t$ is also invariant under rotations and we call $\eta_t$ {\bf isotropic}. Through the paper we make the following assumption on $\QQ$.
\begin{description}
 \item[\textbf{(A1)}] Two independent random subspaces $M,L\in\LG^d_k$ with distribution $\QQ$ are in general position with probability one.
\end{description}
Assumption \textbf{(A1)} is for example fulfilled if $\QQ$ is absolutely continuous with respect to $\nu_k$, see \cite[Theorem 4.4.5]{SW}. We note that under \textbf{(A1)} the flats of $\eta_t$ are almost surely in general position. We also assume henceforth that 
\begin{description}
 \item[\textbf{(A2)}] $1\leq k<d/2$
\end{description}
holds, which ensures that the flats of $\eta_t$ do not intersect each other with probability one (also notice that \textbf{(A2)} implies $d\geq 3$).

Before presenting our main findings in the following three subsections, we introduce some notions and notation used in the present paper. Let us write $\eta_{t,\neq}^2$ for the collection of pairs $(E,F)$ of distinct $k$-flats of $\eta_t$, write $\dist(x,y)$ for the Euclidean distance of two points $x,y\in\RR^d$ and let $\dist(E,F)$ be the \textbf{distance} of two $k$-flats $E,F\in\AG_k^d$, i.e., $\dist(E,F)=\inf\{\dist(x,y):x\in E,\,y\in F\}$. If $E$ and $F$ are in general position, this is the distance of two uniquely determined points $x_E\in E$ and $y_F\in F$ and we call $m(E,F):=(x_E+y_F)/2\in\RR^d$ the \textbf{midpoint} of $E$ and $F$. For two linear subspaces $M,L\in\LG_k^d$ we write $[M,L]$ for the \textbf{subspace determinant} of $M$ and $L$, which is the volume of a parallelepiped generated by two orthonormal bases of $M$ and $L$; cf. \cite[Chapter 14.1]{SW}. Furthermore, we denote in this paper by $V_k(K)$ the intrinsic volume of order $k\in\{0,\ldots,d\}$ of a compact convex set $K\subset\RR^d$; cf. \cite[Chapter 14.2]{SW}. We also write $\kappa_n$ for the volume of the unit ball in $\RR^n$ ($n\geq 1$).

\subsection{The classical proximity}

After these preparations, we can now introduce the \textbf{proximity functional}
$$\pi_t(K,\delta):={1\over 2}\sum_{(E,F)\in\eta_{t,\neq}^2}{\bf 1}\{\dist(E,F)\leq \delta,\, m(E,F)\in K\},$$
where $\delta\in[0,\infty)$ is a fixed threshold, $K$ is a compact and convex subset of $\RR^d$ with $V_d(K)>0$ (called \textbf{convex body} in this paper) and where ${\bf 1}\{\,\cdot\,\}$ is the usual indicator function, which is one if the statement in brackets is fulfilled and zero otherwise. In other words, the functional $\pi_t(K,\delta)$ counts the number of pairs of flats in $\eta_t$ with distance at most $\delta$ and midpoint in $K$. Schneider has calculated in \cite{Schneider99} the mean of $\pi_t(K,\delta)$ for $K$ being the unit ball and $\delta=1$; see also \cite[Theorem 4.4.10]{SW}. More generally, we have the following result.

\begin{theorem}\label{thm:Expectation}
The expectation of $\pi_t(K,\delta)$ is given by
$$\EE \pi_t(K,\delta)=\frac{t^2}{2} \kappa_{d-2k}\delta^{d-2k}\,V_d(K)\int_{\LG^d_k}\int_{\LG^d_k} [M,L]\,\QQ(\dint L)\,\QQ(\dint M).$$
\end{theorem}

\begin{remark}\label{rem:zetakd}\rm
In the isotropic case $\QQ=\nu_k$, we have
\begin{equation}\label{eq:zetakd}
\psi_{d,k}:=\int_{\LG_k^d}\int_{\LG_k^d}[M,L]\,\nu_k(\dint L)\,\nu_k(\dint M)={\kappa_k\kappa_{d-k}\over{d\choose k}\kappa_d},
\end{equation}
which is the content of Corollary 4.5.5 in \cite{Materon}.
\end{remark}

In what follows, we consider a family of increasing observation windows $(K_\varrho)_{\varrho\geq 1}$ with $K_\varrho=\varrho K$ and $K\subset\RR^d$ a convex body and are interested in the asymptotic behaviour of $\pi_t(K_\varrho,\delta)$ as $\varrho\rightarrow\infty$. We first consider the asymptotic variance of $\pi_t(K_\varrho,\delta)$.

\begin{theorem}\label{thm:asymptoticVariance}
It holds that
$$\lim_{\varrho\rightarrow\infty} \frac{\VV\pi_t(K_\varrho,\delta)}{\varrho^{d+k}}=t^3\kappa_{d-2k}^2\delta^{2(d-2k)}\,{\cal I}(K),$$
where $${\cal I}(K)=\int_{\LG^d_k}\int_{M^\perp} V_k\big(K\cap (M+y)\big)^2\,\ell_{M^\perp}(\dint y)\left(\int_{\LG^d_k} [M,L]\,\QQ(\dint L)\right)^2 \QQ(\dint M).$$
\end{theorem}
\begin{remark}\rm
In the case $\QQ=\nu_k$, ${\cal I}(K)$ has an interpretation in terms of the order $k+1$ \textbf{chord-power integral} of $K$, which is defined as $${\cal J}_{k+1}(K):=\int_{\AG_1^d}V_1(K\cap G)^{k+1}\,\mu_1(\dint G),$$ where $\mu_1$ is the Haar measure on $\AG_1^d$ normalized as in \cite{SW}. Indeed, we first notice that the rotational average $\int_{\LG^d_k} [M,L]\,\nu_k(\dint L)$ does not depend on $M$; cf. \cite[Corollary 4.5.5]{Materon}. Then identity (8.57) in \cite{SW} implies that $${\cal I}(K)={\kappa_k\over k+1}\psi_{d,k}^2\,{\cal J}_{k+1}(K)$$ with $\psi_{d,k}$ as in \eqref{eq:zetakd}.
\end{remark}

Having investigated the expectation and the asymptotic variance of the proximity functional $\pi_t(K,\delta)$, we turn now to the central limit problem. Let the family $(K_\varrho)_{\varrho\geq 1}$ of convex bodies be as above.

\begin{theorem}\label{thm:CLT}
Let $\cN$ be a standard Gaussian random variable. Then there is a constant $C$ depending on $K$, $\delta$ and $t$ such that
$$\sup_{x\in\RR}\left|\PP\left(\frac{\pi_t(K_\varrho,\delta)-\EE\pi_t(K_\varrho,\delta)}{\sqrt{\VV\pi_t(K_\varrho,\delta)}}\leq x\right)-\PP(\cN\leq x)\right|\leq C\,\varrho^{-{d-k\over 2}}$$ for $\varrho\geq 1$. In particular, we have the convergence in distribution $$\frac{\pi_t(K_\varrho,\delta)-\EE\pi_t(K_\varrho,\delta)}{\sqrt{\VV\pi_t(K_\varrho,\delta)}}\overset{d}{\longrightarrow}\cN\quad{\rm as}\quad \varrho\rightarrow\infty.$$
\end{theorem}

\subsection{Small distances}

In the previous theorems, we have considered the number of midpoints of pairs of flats in a sequence of increasing observation windows, which have distance below a given threshold $\delta$. A further natural question is to ask for the shortest or, more generally, the $m$-th shortest distance between two flats. To present the result, let $(K_\varrho)_{\varrho\geq 1}$ be a family of convex bodies as above. We denote by
\begin{equation}\label{eq:PPXI}
\xi_\varrho^{(K,t)} = \{\dist(E,F): (E,F)\in\eta^2_{t,\neq} \text{ and } m(E,F)\in K_\varrho\}
\end{equation}
the set of all distances between pairs of flats having a midpoint in $K_\varrho$ (we count each value $\dist(E,F)$ only once, although $(E,F)$ and $(F,E)$ are both elements of $\eta^2_{t,\neq}$). Formally, $\xi_\varrho^{(K,t)}$ can be considered as a point process on the positive real half-line $\RR_+$. By $D_m^{(K_\varrho,t)}$ we denote the $m$-th smallest element in $\xi_\varrho^{(K,t)}$ according to the natural ordering on $\RR_+$. The following theorem describes the asymptotic distributions of $D_m^{(K_\varrho,t)}$ and $\xi_\varrho^{(K,t)}$ as the window size tends to infinity.

\begin{theorem}\label{thm:PointProcess}
Define
\begin{equation*}
\beta=\frac{t^2}{2}\kappa_{d-2k}V_d(K) \int_{\LG^d_k}\int_{\LG^d_k} [M,L]\,\QQ(\dint L)\,\QQ(\dint M).
\end{equation*}
For every $u\geq0$, there exists a constant $C_u$ also depending on $K$ and $t$ such that
$$\left|\PP\big(\varrho^{d/(d-2k)}D_m^{(K_\varrho,t)}>u\big)-e^{-\beta u^{(d-2k)}}\sum_{i=0}^{m-1}\frac{(\beta u^{(d-2k)})^i}{i!}\right|\leq C_u\,\varrho^{-\frac{d-k}{2}}$$
for $m\in\{1,2,3,\ldots\}$ and $\varrho\in[1,\infty)$. Moreover, the family $\big(\varrho^{d/(d-2k)}\xi_\varrho^{(K,t)}\big)_{\varrho\geq 1}$ of rescaled point processes converges in distribution to a Poisson point process on $\RR_+$ with the intensity measure
$$\nu(A)=\beta(d-2k)\int_A u^{d-2k-1}\,\dint u,\quad A\subset\RR_+\ \text{ a Borel set}.$$
\end{theorem}

\begin{remark}\rm
We notice that $$\beta={t^2\over 2}\kappa_{d-2k}\,\psi_{d,k}\,V_d(K)$$ with $\psi_{d,k}$ given by \eqref{eq:zetakd} in the case where $\QQ=\nu_k$ is the invariant probability measure on $\LG_k^d$.
\end{remark}

\begin{remark}\rm
In \cite{SchulteThaele12} a similar problem was considered. Namely, for a pair $(E,F)$ of flats of a Poisson $k$-flat process with $\QQ=\nu_k$ hitting a convex body $K$, the distance $\dist_K(E,F)$ was defined as $$\dist_K(E,F):=\min_{x\in E\cap K,y\in F\cap K}\dist(x,y)$$ and it was shown that for increasing intensity the ordered distances converge to an inhomogeneous Poisson point process similar to that in Theorem \ref{thm:PointProcess}. The fact that increasing the intensity is up to a factor the same as increasing the window size implies that that the normalization $\varrho^{d/(d-2k)}$ in Theorem \ref{thm:PointProcess} is the same as in \cite{SchulteThaele12}. The constants $\beta$, however, are different in both settings since different pairs of flats and different approaches to measure the distance between two flats are considered.
\end{remark}

\subsection{Distances around a positive value}

The previous result describes the behaviour of very small distances and it is natural also to consider large distances. However, the maximal distance (and thus also the $m$-th maximal distance for any $m\in\{1,2,3,\ldots\}$) of two flats having their midpoint in a test set $K$ is not well defined since
\begin{equation}\label{eq:Maximum}
\sup_{\stackrel{(E,F)\in\eta_{t,\neq}^2}{m(E,F)\in K}} \dist(E,F)=\infty\quad \text{almost surely};
\end{equation}
see Section \ref{sec:Extremes} for a proof.

To overcome this difficulty and in order to complete the picture, we fix some $\sigma>0$ and consider the asymptotic behaviour of the point process $\xi_\varrho^{(K,t)}$ defined by \eqref{eq:PPXI} around $\sigma$. By $\overline{D}^{(K_\varrho,t,\sigma)}_m$ and $\underline{D}_m^{(K_\varrho,t,\sigma)}$, $m\in\{1,2,3,\ldots\}$, we denote the $m$-th element of $\xi_\varrho^{(K,t)}$ greater or less than $\sigma$, respectively. 

\begin{theorem}\label{thm:Maximum}
Let $\sigma>0$ and define
\begin{equation}\label{eq:betaMAX}
\beta=\frac{t^2}{2} (d-2k)\kappa_{d-2k}\,\sigma^{d-2k-1}\,V_d(K)\int_{\LG^d_k}\int_{\LG^d_k} [M,L]\,\QQ(\dint L)\,\QQ(\dint M).
\end{equation}
For every $u\geq0$, there is a constant $C_u$ also depending on $K$, $t$ and $\sigma$ such that
$$\left|\PP\big(\varrho^d\,(\overline{D}_m^{(K_\varrho,t,\sigma)}-\sigma)>u\big)-e^{-\beta u}\sum_{i=0}^{m-1} \frac{(\beta u)^i}{i!}\right|\leq C_u\,\varrho^{-\frac{d-k}{2}}$$ and
$$\left|\PP\big(-\varrho^d\,(\underline{D}_m^{(K_\varrho,t,\sigma)}-\sigma)>u\big)-e^{-\beta u}\sum_{i=0}^{m-1} \frac{(\beta u)^i}{i!}\right|\leq C_u\,\varrho^{-\frac{d-k}{2}}$$
for $m\in\{1,2,3,\ldots\}$ and $\varrho\in[1,\infty)$. Moreover, the family $\big(\varrho^d\,(\xi_\varrho^{(K,t)}-\sigma)\big)_{\varrho\geq 1}$ of rescaled and shifted point processes converges in distribution to a homogeneous Poisson point process on $\RR$ with intensity $\beta$.
\end{theorem}

\begin{remark}\rm
In the case where $\QQ$ is the invariant probability measure $\nu_k$ on $\LG_k^d$ we have that $$\beta={t^2\over 2}(d-2k)\kappa_{d-2k}\,\sigma^{d-2k-1}\,\psi_{d,k}\,V_d(K)$$ with $\psi_{d,k}$ given by \eqref{eq:zetakd}.
\end{remark}

Theorem \ref{thm:PointProcess} and Theorem \ref{thm:Maximum} show the remarkable fact that very small distances near zero behave quite different compared with the distances around (i.e., above or below) every positive value $\sigma$. Indeed, in Theorem \ref{thm:PointProcess} an inhomogeneous Poisson point process on $\RR_+$ appears after normalization with $\varrho^{d/(d-2k)}$, whereas in Theorem \ref{thm:Maximum} a homogeneous Poisson point process on the whole real line shows up in the limit after rescaling with $\varrho^d$ and the latter can be interpreted as the superposition of two independent homogeneous Poisson point process on $\RR_+$ (for the distances greater than $\sigma$) and on $\RR_-$ (for the distances less than $\sigma$).

\section{Background material on chaos decompositions}\label{sec:chaos}

We let $\eta_t$ be a Poisson point process on $\AG_k^d$ with intensity measure $\Theta_t$ given by \eqref{eq:intensitymeasure} and assume that \textbf{(A1)} and \textbf{(A2)} are satisfied. Given $n\in\{1,2\}$ we write $L^2(\Theta_t^n)$ for the collection of functions $f:(\AG_k^d)^n\rightarrow\RR$ such that $$||f||_n:=\Big(\int_{(\AG_k^d)^n}f^2\,\dint\Theta_t^n\Big)^{1/2}<\infty$$ and $L_{\rm sym}^2(\Theta_t^2)$ for the subspace of $L^2(\Theta_t^2)$ consisting of functions that are invariant under permutation of the two arguments, so called symmetric functions (formally, we also have $L_{\rm sym}^2(\Theta_t)=L^2(\Theta_t)$).

For $f\in L_{\rm sym}^2(\Theta_t^n)$ we let $I_n(f)$ be the (multiple) \textbf{Wiener-It\^o integral} of $f$ with respect to the \textbf{compensated Poisson process} $\hat\eta_t:=\eta_t-\Theta_t$ (to make sense of the definition of $\hat\eta_t$, $\eta_t$ has to be interpreted here as a random point measure so that the difference $\hat\eta_t-\Theta_t$ is well defined). This is to say, $$I_1(f)=\int_{\AG_k^d}f(E)\,\hat\eta_t(\dint E)\qquad{\rm if}\ n=1$$ and $$I_2(f)=\int_{\AG_k^d}\int_{\AG_k^d}f(E,F)\,\hat\eta_t(\dint F)\,\hat\eta_t(\dint E)\qquad{\rm if}\ n=2;$$ cf. \cite{LP,PeccatiTaqqu08,PTSU10}. These stochastic integrals satisfy the following properties: it holds that
\begin{equation}\label{eq:rule1}
\EE I_n(f)=0\quad{\rm and}\quad\EE I_n(f)^2= n!||f||_n^2
\end{equation}
for $f\in L^2_{\rm sym}(\Theta_t^n)$ and
\begin{equation}\label{eq:rule2}
\EE[I_1(f_1)I_2(f_2)]=0
\end{equation}
for $f_1\in L^2(\Theta_t)$ and $f_2\in L_{\rm sym}^2(\Theta_t^2)$.

Let $g:(\AG_k^d)^2\rightarrow\RR$ be integrable with respect to $\Theta_t^2$ and be invariant under permutation of its two arguments.
We define $$U:={1\over 2}\sum_{(E,F)\in\eta_{t,\neq}^2}g(E,F)$$ and assume that $U$ is square integrable with respect to the distribution of $\eta_t$. In this case, the random variable $U$ is a so-called \textbf{Poisson U-statistic of order two}.
It is a crucial fact that $U$ can be written as
\begin{equation}\label{eq:chaos}
U = \EE U + I_1(f_1)+I_2(f_2) 
\end{equation}
with
\begin{equation}\label{eq:campbell}
\EE U = {1\over 2}\int_{(\AG_k^d)^2}g(E,F)\,\Theta_t^2\big(\dint(E,F)\big)
\end{equation}
by the classical Slivnyak-Mecke formula \cite[Theorem 3.2.5]{SW} and with $f_1\in L_{\rm sym}^2(\Theta_t)$ and $f_2\in L_{\rm sym}^2(\Theta_t^2)$ given by
\begin{eqnarray*}
f_1(E) &=& \int_{\AG_k^d}g(E,F)\,\Theta_t(\dint F),\\
f_2(E,F) &=& {1\over 2}\,g(E,F);
\end{eqnarray*}
cf. Lemma 3.5 in \cite{RS12}. The representation \eqref{eq:chaos} is called the \textbf{Wiener-It\^o chaos decomposition} of $U$ and we call $f_1$ and $f_2$ its \textbf{kernels}. This decomposition is a very powerful tool, which will be used extensively in our proofs below. In particular, squaring \eqref{eq:chaos} and using the computation rules \eqref{eq:rule1} and \eqref{eq:rule2}, we find the \textbf{variance formula}
\begin{equation}\label{eq:variance}
\VV U = ||f_1||_1^2+2||f_2||_2^2.
\end{equation}
This will be essential in the proof of Theorem \ref{thm:asymptoticVariance}.

For more details on Poisson U-statistics (of arbitrary order) we refer to \cite{LRPeccati1,LRPeccati2,LPST,RS12}. Poisson functionals that are a sum of a first and a second order Wiener-It\^o integral, such as the functional $U$ above, were also investigated in \cite{PeccatiTaqqu08}.

\section{Proof of Theorems \ref{thm:Expectation} and \ref{thm:asymptoticVariance}}

\subsection{A preparatory lemma}

In order to simplify our notation, we define for $E,F\in\AG^d_k$, $$h(E,F)={\bf 1}\{m(E,F)\in K,\,\dist(E,F)\leq \delta\},$$ where $K$ is a convex body and where $\delta>0$. 

\begin{lemma}\label{lem:Integral}
Let $M,L\in\LG^d_k$ be in general position and define $W=M+L$. Then
\begin{equation}\label{eq:eqlemma}
\begin{split}
\int_{L^\perp} h(M,L+x)\, &\ell_{L^\perp}(\dint x)\\ =[M,L]\,&\int_{W^\perp} {\bf 1}\{\|x\|\leq\delta\}\,V_k\big((K-(x/2))\cap M\big)\, \ell_{W^\perp}(\dint x).
\end{split}
\end{equation}
\end{lemma}

\paragraph{Proof:} By decomposing $x\in L^\perp$ in $x=x_1+x_2$ with $x_1\in L^\perp\cap W^\perp=W^\perp$ and $x_2\in L^\perp\cap W$, we obtain
\begin{equation}\label{eq:lemmaproof1}
\begin{split}
\int_{L^\perp} h(M,L+x)\, &\ell_{L^\perp}(\dint x)\\ = \int_{W^\perp}\int_{L^\perp\cap W} &h(M,L+x_1+x_2)\,\ell_{L^\perp\cap W}(\dint x_2)\,\ell_{W^\perp}(\dint x_1). 
\end{split} 
\end{equation}
By the definition of $x_1$ and $x_2$, $M$ and $L+x_2$ intersect in a unique point $z\in W$, $m(M,L+x_1+x_2)=z+(x_1/2)$ and $\dist(M,L+x_1+x_2)=\|x_1\|$.

Let $B_M$ and $B_{L^\perp\cap W}$ be matrices whose columns form orthonormal bases of $M$ and $L^\perp\cap W$, respectively. Rewrite $x_2$ as $x_2=B_{L^\perp\cap W}\tilde{x}$ with $\tilde{x}\in\RR^k$ and replace integration over $L^\perp\cap W$ in \eqref{eq:lemmaproof1} by integration over $\RR^k$. Moreover, we notice that the intersection point $z$ of $M$ and $L+x_2$ has the representation $z=B_M\tilde{z}$, where $\tilde{z}\in\RR^k$ is the solution of
$$B^T_{L^\perp\cap W}B_M \tilde{z}=B^T_{L^\perp\cap W}x_2=B^T_{L^\perp\cap W}B_{L^\perp\cap W}\tilde{x}=\tilde{x}.$$ This implies that 
\begin{equation}\label{eq:lemmaproof2}
\tilde{z}=\big(B^T_{L^\perp\cap W}B_M\big)^{-1}\tilde{x} .
\end{equation}
Using the representation of $x_2$, we now write the inner integral in \eqref{eq:lemmaproof1} as
\begin{equation}\label{eq:exproof1}
\begin{split}
& \int_{L^\perp\cap W} h(M,L+x_1+x_2)\,\ell_{L^\perp\cap W}(\dint x_2)\\
& = {\bf 1}\{\|x_1\|\leq \delta\}\,\int_{\RR^k} {\bf 1}\{m(M,L+x_1+B_M\tilde{x})\in K\}\, \dint\tilde{x}.
\end{split}
\end{equation}
Continuing by using \eqref{eq:lemmaproof2}, we find
\begin{equation}\label{eq:exproof2}
\begin{split}
& \int_{\RR^k} {\bf 1}\{m(M,L+x_1+B_M\tilde{x})\in K\}\, \dint\tilde{x}\\ 
& = \int_{\RR^k}{\bf 1}\{B_M(B^T_{L^\perp\cap W}B_M)^{-1}\tilde{x}\in\big(K-(x_1/2)\big)\cap M\}\,\dint\tilde{x}\\
& = \int_{\RR^k} {\bf 1}\{\tilde{x}\in B^T_{L^\perp\cap W}B_M B_M^{-1}(K-(x_1/2))\cap M\}\,\dint\tilde{x}\\
& = V_k\big(B^T_{L^\perp\cap W}B_M B_M^{-1}(K-(x_1/2))\cap M\big).
\end{split}
\end{equation}
Combining \eqref{eq:exproof1} with \eqref{eq:exproof2} and using the fact that $\det(B^T_{L^\perp\cap W}B_M)=[M,L]$, we arrive at
$$\hspace{-2cm}\int_{L^\perp\cap W} h(M,L+x_1+x_2)\,\ell_{L^\perp\cap W}(\dint x_2)$$ $$\hspace{2cm}=[M,L]\,{\bf 1}\{\|x_1\|\leq \delta\}V_k\big((K-(x_1/2))\cap M\big).$$ Integration with respect to $W^\perp$ finally yields \eqref{eq:eqlemma}.
\hfill $\Box$

\subsection{Proof of Theorem \ref{thm:Expectation}}

By \eqref{eq:campbell}, the expectation of $\pi_t(K,\delta)$ is given by 
$$\EE \pi_t(K,\delta)={1\over 2}\int_{\AG_k^d}\int_{\AG_k^d}h(E,F)\,\Theta_t(\dint F)\,\Theta_t(\dint E).$$ A glance at \eqref{eq:intensitymeasure} shows that this equals
\begin{equation}\label{eq:thm1proof1}
\frac{t^2}{2}\int_{\LG^d_k}\int_{\LG^d_k}\int_{M^\perp}\int_{L^\perp} h(M+y,L+x)\,\ell_{L^\perp}(\dint x)\,\ell_{M^\perp}(\dint y)\,\QQ(\dint L)\,\QQ(\dint M).
\end{equation}
We evaluate now the inner double integral in \eqref{eq:thm1proof1}. Using Lemma \ref{lem:Integral}, we obtain
\begin{equation*}
\begin{split}
& \int_{M^\perp}\int_{L^\perp} h(M+y,L+x)\,\ell_{L^\perp}(\dint x)\,\ell_{M^\perp}(\dint y)\\
& = \int_{M^\perp}\int_{L^\perp} {\bf 1}\{m(M,L+x)\in K-y,\,\dist(M,L+x)\leq \delta\}\, \ell_{L^\perp}(\dint x)\,\ell_{M^\perp}(\dint y)\\
& = [M,L]\int_{M^\perp}\int_{W^\perp} V_k\big((K-y-(x/2))\cap M\big)\,{\bf 1}\{\|x\|\leq\delta\}\,\ell_{W^\perp}(\dint x)\,\ell_{M^\perp}(\dint y),
\end{split}
\end{equation*}
where $W=L+M$. Fubini's theorem further implies that
\begin{equation*}
\begin{split}
& [M,L]\int_{M^\perp}\int_{W^\perp} V_k\big((K-y-(x/2))\cap M\big)\,{\bf 1}\{\|x\|\leq\delta\}\,\ell_{W^\perp}(\dint x)\,\ell_{M^\perp}(\dint y)\\
& = [M,L]\int_{W^\perp}\int_{M^\perp} V_k\big((K-y-(x/2))\cap M\big)\,{\bf 1}\{\|x\|\leq\delta\}\,\ell_{M^\perp}(\dint y)\,\ell_{W^\perp}(\dint x)\\
& = [M,L]\,V_d(K)\,\kappa_{d-2k}\delta^{d-2k},
\end{split}
\end{equation*}
which in view of \eqref{eq:thm1proof1} completes the proof. \hfill $\Box$\\

\begin{remark}\rm
The above proofs of Lemma \ref{lem:Integral} and Theorem \ref{thm:Expectation} are very similar to that of the main result in \cite{Schneider99} and use the same ideas, generalized to a slightly more general setting. We decided to state the first part as lemma since \eqref{eq:eqlemma} is applied several times below.
\end{remark}

\subsection{Proof of Theorem \ref{thm:asymptoticVariance}}

First, \eqref{eq:chaos} implies that the proximity functional $\pi_t(K,\delta)$ has chaos decomposition
$$\pi_t(K,\delta)=\EE\pi_t(K,\delta)+I_1(f_1^{(K,\delta,t)})+I_2(f_2^{(K,\delta,t)}).$$ Here, the kernels $f_n^{(K\delta,t)}$ $(n=1,2)$ are given by
\begin{equation}\label{eq:f1}
f_1^{(K,\delta,t)}(M+y) = \int_{\AG_k^d} h(M+y,F)\,\Theta_t(\dint F)
\end{equation}
for $M\in\LG^d_k$ and $y\in M^\perp$ and
\begin{equation}\label{eq:f2}
f_2^{(K,\delta,t)}(E,F) = \frac{1}{2}\,{\bf 1}\{m(E,F)\in K,\,\dist(E,F)\leq \delta\}
\end{equation}
for $E,F\in\AG^d_k$, respectively. Now, the variance formula \eqref{eq:variance} implies that
\begin{equation}\label{eq:VariancePi}
\VV \pi_t(K,\delta)=\|f_1^{(K,\delta,t)}\|_1^2+2\|f_2^{(K,\delta,t)}\|_2^2.
\end{equation}
We determine the asymptotic behaviour of the right hand side in \eqref{eq:VariancePi}. For the second term we find
\begin{eqnarray*}
\|f_2^{(K,\delta,t)}\|_2^2 & = & \frac{1}{4}\int_{\AG^d_k}\int_{\AG^d_k} {\bf 1}\{m(E,F)\in K,\, \dist(E,F)\leq\delta\}^2\, \Theta_t(\dint E)\,\Theta_t(\dint F)\\
& = &\frac{t^2}{4}\,\kappa_{d-2k}\,\delta^{d-2k}\,V_d(K)\int_{\LG^d_k}\int_{\LG_k^d} [L,M]\,\QQ(\dint L)\,\QQ(\dint M)
\end{eqnarray*}
by using the formula for $\EE\pi_t(K,\delta)$ in Theorem \ref{thm:Expectation}. Thus,
\begin{equation}\label{eq:proofvariancesecondterm}
\begin{split}
& \lim_{\varrho\rightarrow\infty}\frac{\|f^{(K_\varrho,\delta,t)}_2\|_2^2}{\varrho^{d+k}}\\
& =\lim_{\varrho\rightarrow\infty}\frac{t^{2}}{4}\,\kappa_{d-2k}\,\delta^{d-2k}V_d(K)\varrho^{-k}\int_{\LG^d_k}\int_{\LG^d_k}[M,L]\,\QQ(\dint L)\,\QQ(\dint M)\\
& =0.
\end{split} 
\end{equation}
We continue with the first term in \eqref{eq:VariancePi} and observe that \eqref{eq:intensitymeasure} and Lemma \ref{lem:Integral} imply that
\begin{equation*}
\begin{split}
& f_1^{(K,\delta,t)}(M+y) = \int_{\AG_k^d} h(M+y,F)\,\Theta_t(\dint F)\\
& = t\int_{\LG_k^d}\int_{L^\perp} h(M+y,L+x)\,\ell_{L^\perp}(\dint x)\,\QQ(\dint L)\\
& = t \int_{\LG^d_k} \int_{W^\perp} V_k\big((K-y-(x/2))\cap M\big)\,{\bf 1}\{\|x\|\leq \delta\}\,\ell_{W^\perp}(\dint x)\,[M,L]\,\QQ(\dint L),
\end{split}
\end{equation*}
where, as before, $W=L+M$. We now observe that the scaling relation
\begin{equation}\label{eq:proofvariancescaling}
f_1^{(K_\varrho,\delta,t)}(M+y)= \varrho^{d-k}\,f_1^{(K,\delta/\varrho,t)}(M+(y/\varrho))
\end{equation}
holds. Indeed, from a simple change of variables and from the fact that ${\rm dim}W^\perp=d-2k$ it follows that
\begin{equation*}
\begin{split}
& \int_{W^\perp} V_k\big((K_\varrho-y-(x/2))\cap M\big)\,{\bf 1}\{\|x\|\leq \delta\}\,\ell_{W^\perp}(\dint x)\\
& = \varrho^k \int_{W^\perp}  V_k\big((K-(y/\varrho)-(x/2\varrho))\cap M\big)\,{\bf 1}\{\|x\|\leq \delta\}\,\ell_{W^\perp}(\dint x)\\
& = \varrho^{d-k}\int_{W^\perp}  V_k\big((K-(y/\varrho)-(x/2))\cap M\big)\,{\bf 1}\{\|x\|\leq \delta/\varrho\}\,\ell_{W^\perp}(\dint x),
\end{split}
\end{equation*}
which shows \eqref{eq:proofvariancescaling}. As a consequence, we have
\begin{equation}\label{eq:f1asymptotic}
\begin{split}
& \|f_1^{(K_\varrho,\delta,t)}\|^2_1\\
& = \varrho^{2(d-k)}\,t\int_{\LG^d_k} \int_{M^\perp} f_1^{(K,\delta/\varrho,t)}(M+(y/\varrho))^2\,\ell_{M^\perp}(\dint y)\,\QQ(\dint M)\\
& = \varrho^{3(d-k)}\,t\int_{\LG^d_k}\int_{M^\perp} f_1^{(K,\delta/\varrho,t)}(M+y)^2\,\ell_{M^\perp}(\dint y)\,\QQ(\dint M).
\end{split}
\end{equation}
Moreover, the dominated convergence theorem implies that
\begin{equation*}
\begin{split}
& \lim_{\varrho\rightarrow\infty} \varrho^{d-2k} f_1^{(K,\delta/\varrho,t)}(M+y)\\
& = t \int_{\LG^d_k} \lim_{\varrho\rightarrow\infty}\varrho^{d-2k}\int_{W^\perp} V_k\big((K-y-(x/2))\cap M\big)\\
& \qquad\qquad\qquad\qquad\times\,{\bf 1}\{\|x\|\leq \delta/\varrho\}\,\ell_{W^\perp}(\dint x)\,[M,L]\,\QQ(\dint L)\\
& = t\,\kappa_{d-2k}\,\delta^{d-2k}\,V_k\big((K-y)\cap M\big) \int_{\LG^d_k} [M,L]\,\QQ(\dint L)\\
& = t\,\kappa_{d-2k}\,\delta^{d-2k}\,V_k\big(K\cap (M+y)\big) \int_{\LG^d_k} [M,L]\,\QQ(\dint L).
\end{split}
\end{equation*}
Combining this with \eqref{eq:f1asymptotic}, writing $\varrho^{3(d-k)}$ as $\varrho^{2(d-2k)}\varrho^{d+k}$ and applying the dominated convergence Theorem once again, yields
$$\lim_{\varrho\rightarrow\infty} \frac{\|f_1^{(K_\varrho,\delta,t)}\|_1^2}{\varrho^{d+k}}=t^3\,\kappa_{d-2k}^2\,\delta^{2(d-2k)}\,{\cal I}(K)$$ with ${\cal I}(K)$ as in the statement of the theorem. This together with the asymptotic behaviour \eqref{eq:proofvariancesecondterm} of the second term in the variance expansion \eqref{eq:VariancePi} of the proximity functional proves the claim. \hfill $\Box$

\section{Proof of Theorem \ref{thm:CLT}}

\subsection{A general bound}

For two random variables $Y$ and $Z$ define the Kolmogorov distance $d_K(Y,Z)$ by
$$d_K(Y,Z)=\sup_{x\in\RR}|\PP(Y\leq x)-\PP(Z\leq x)|.$$ This is to say, $d_K(Y,Z)$ is the supremum norm of the difference between the distribution functions of $Y$ and $Z$. We consider a second-order Poisson U-statistic
$$U=\frac{1}{2}\sum_{(E,F)\in\eta^2_{t,\neq}} g(E,F),$$
where we assume that $g$ is bounded, symmetric and satisfies
$$\Theta_t^2\big(\{(E,F)\in\AG^d_k\times\AG^d_k: g(E,F)\neq 0\}\big)<\infty,\qquad t>0.$$
We denote the kernels of the chaos decomposition of $U$ given in \eqref{eq:chaos} by $f_1$ and $f_2$ and define $M_{11}$ by
$$M_{11} = \int_{\AG^d_k} f_1(E)^4\,\Theta_t(\dint E).$$
We also define $M_{12}$ by
\begin{eqnarray*}
M_{12} & = & 8 \int_{(\AG^d_k)^3} f_1(E_1) f_2(E_1,E_2) f_1(E_3) f_2(E_2,E_3)\, \Theta_t^3\big(\dint (E_1,E_2,E_3)\big)\\
&& +\,4\int_{(\AG^d_k)^2} f_1(E_1) f_2(E_1,E_2) f_1(E_2) f_2(E_1,E_2)\, \Theta_t^2\big(\dint (E_1,E_2)\big)
\end{eqnarray*}
and finally $M_{22}$ by
\begin{eqnarray*}
M_{22} & = & 48 \int_{(\AG^d_k)^4} f_2(E_1,E_2) f_2(E_2,E_3) f_2(E_3,E_4)\\
& & \qquad\qquad\qquad\qquad\qquad\times\,f_2(E_4,E_1)\,\Theta_t^4\big(\dint (E_1,E_2,E_3,E_4)\big)\\
&& \qquad+\, 96 \int_{(\AG^d_k)^3} f_2(E_1,E_2) f_2(E_1,E_3) f_2(E_1,E_3) \\
& & \qquad\qquad\qquad\qquad\qquad\times\,f_2(E_2,E_3)\,\Theta_t^3\big(\dint (E_1,E_2,E_3)\big)\\
&& \qquad\qquad+\,8\,\int_{(\AG^d_k)^2}f_2(E_1,E_2)^4\,\Theta_t^2\big(\dint (E_1,E_2)\big).
\end{eqnarray*}
We can now rephrase a special situation of Theorem 4.2 in \cite{Schulte2012}.
 
\begin{proposition}\label{prop:clt}
Let $\cN$ be a standard Gaussian random variable. Then 
\begin{equation}\label{eq:Kolmogorov}
d_K\left(\frac{U-\EE U}{\sqrt{\VV U}},\cN\right) \leq 1088\, \frac{\sqrt{M_{11}}+\sqrt{M_{12}}+\sqrt{M_{22}}}{\VV U}.
\end{equation}
\end{proposition}

\subsection{Proof of Theorem \ref{thm:CLT}}

Let us introduce the abbreviation
$$\pi_\varrho:=\pi_t(K_\varrho,\delta)=\frac{1}{2}\sum_{(E,F)\in\eta^2_{t,\neq}}{\bf 1}\{m(E,F)\in K_\varrho,\,\dist(E,F)\leq \delta\}.$$
Since $t$, $K$ and $\delta$ are fixed in the following, we suppress this dependency in our notation. We further let $f^{(\varrho)}_1$ and $f^{(\varrho)}_2$ be the kernels of the Wiener-It\^o chaos decomposition of $\pi_\varrho$ given by \eqref{eq:f1} and \eqref{eq:f2}, respectively. In the following, we prove Theorem \ref{thm:CLT} by bounding the right hand side of \eqref{eq:Kolmogorov} for the Poisson U-statistic $U=\pi_\varrho$. 

\paragraph{Step 1: Two inequalities for $f_1^{(\varrho)}$ and $f_2^{(\varrho)}$.} We show that 
\begin{equation}\label{eq:proofCLT1}
f_1^{(\varrho)}(E)\leq C\varrho^k\ \text{ and } \ \Theta_t\big(\{F\in\AG^d_k: f_2^{(\varrho)}(E,F)\neq 0\}\big)\leq C\varrho^k
\end{equation}
for all $E\in\AG^d_k$, where $C=t\kappa_k\kappa_{d-2k}\delta^{d-2k}\,(\diam(K)/2)^k$. The first bound is a consequence of \eqref{eq:f1}, Lemma \ref{lem:Integral} and the inequality $V_k(\tilde{K})\leq\kappa_k(\diam(\tilde{K})/2)^k$ from \cite[page 76]{BF} for a convex body $\tilde{K}\subset\RR^k$. Indeed, writing $E=M+y$ and $W=M+L$, it holds that
\begin{equation*}
\begin{split}
& f_1^{(\varrho)}(E)=f_1^{(\varrho)}(M+x)\\
&=t\int_{\LG_k^d}[M,L]\int_{W^\perp}{\bf 1}\{\|y\|\leq\delta\}\,V_k\big((K_\varrho-(y/2))\cap M\big)\,\ell_{W^\perp}(\dint y)\,\QQ(\dint L)\\
&\leq t\int_{\LG_k^d}[M,L]\,\int_{W^\perp}{\bf 1}\{\|y\|\leq\delta\}\,\kappa_k\,({\rm diam}(K_\varrho)/2)^k\,\ell_{W^\perp}(\dint y)\,\QQ(\dint L)\\
&\leq t\kappa_k\kappa_{d-2k}\delta^{d-2k}\,({\rm diam}(K)/2)^k\,\varrho^k,
\end{split}
\end{equation*}
where we have used additionally the fact that $[M,L]\leq 1$. To show the second inequality, we apply once more Lemma \ref{lem:Integral} to see that
\begin{equation*}
\begin{split}
& \Theta_t\big(\{F\in\AG^d_k: f_2^{(\varrho)}(E,F)\neq 0\}\big)\\
& = t\int_{\LG^d_k}\int_{W^\perp} {\bf 1}\{\|x\|\leq\delta\}\,V_k\big((K-y-(x/2))\cap M\big)\,\ell_{W^\perp}(\dint x) [M,L]\,\QQ(\dint L). 
\end{split}
\end{equation*}
Now, the same argumentation as above yields the second part of \eqref{eq:proofCLT1}.

\paragraph{Step 2: Completing the proof.}
Let $B^d_\delta$ be the $d$-dimensional centred ball with radius $\delta$ and denote by $+$ the usual Minkowski sum.

All integrands occurring in $M_{11}$, $M_{12}$ and $M_{22}$ have the structure that after choosing the first $k$-flat $E$ hitting $K_\varrho+ B_\delta^d$, the second flat must be in the set $\{F\in\AG^d_k: f_2^{(\varrho)}(E,F)\neq 0\}$ or the integrand is zero otherwise. For the remaining flats there are similar conditions so that, by Step 1, the measure of the support of each integrand is at most
$$\Theta_t\big([K_\varrho+ B^d_\delta]\big)\,(C\varrho^k)^{m-1}\leq\Theta_t\big([K+ B^d_\delta]\big)\,C^{m-1}\varrho^{d+(m-2)k}$$
for $\varrho\geq 1$. Here, $m\in\{1,2,3,4\}$ is the number of $k$-flats the integration runs over and for a set $A\subset\RR^d$, $[A]$ stands for the collection of $k$-flats that have non-empty intersection with $A$. Combining this with the fact that $f_1^{(\varrho)}\leq C \varrho^{k}$, recall \eqref{eq:proofCLT1}, and $f_2^{(\varrho)}\leq \frac{1}{2}$, we obtain
\begin{eqnarray*}
M_{11} & \leq & C^4\,\Theta_t\big([K+ B^d_\delta]\big)\,\varrho^{d+3k},\\
M_{12} & \leq & 2C^4\,\Theta_t\big([K+ B^d_\delta]\big)\,\varrho^{d+3k}+C^3\,\Theta_t\big([K+ B^d_\delta]\big)\,\varrho^{d+2k},\\
M_{22} & \leq & 3C^3\,\Theta_t\big([K+ B^d_\delta]\big)\,\varrho^{d+2k}+ 6C^2\,\Theta_t\big([K+ B^d_\delta]\big)\,\varrho^{d+k}\\
& & \qquad\qquad\qquad + \frac{C}{2}\,\Theta_t\big([K+B^d_\delta]\big)\,\varrho^d.  
\end{eqnarray*}
On the other hand, Theorem \ref{thm:asymptoticVariance} tells us that $\VV \pi_\varrho$ is asymptotically of order $\varrho^{d+k}$, so that $\sqrt{M_{11}}/\VV \pi_\varrho$, $\sqrt{M_{12}}/\VV \pi_\varrho$ and $\sqrt{M_{22}}/\VV \pi_\varrho$ are of order $\varrho^{-(d-k)/2}$ or less and Proposition \ref{prop:clt} implies Theorem \ref{thm:CLT}.\hfill $\Box$

\section{Proofs of Theorem \ref{thm:PointProcess}, Theorem \ref{thm:Maximum} and Equation (\ref{eq:Maximum})}\label{sec:Extremes}

\subsection{An auxiliary limit theorem}

We consider the following general setting. Let $(g_\varrho)_{\varrho\geq 1}$ be a family of symmetric functions $g_\varrho: (\AG^d_k)^2\rightarrow\RR$ satisfying $\Theta_t^2\big(g_\varrho^{-1}([-u,u])\big)<\infty$ for all $u\geq 0$ (this will always be the case in our applications below). Next, we define a point process
$$\xi_\varrho=\{g_\varrho(E,F): (E,F)\in\eta^2_{t,\neq},\, m(E,F)\in K_\varrho\}$$
on $\RR$, where we count the point $g_\varrho(E,F)=g_\varrho(F,E)$ only once and where the family $(K_\varrho)_{\varrho\geq 1}$ is a family of convex bodies as in Section \ref{sec:results} (that $\xi_\varrho$ is indeed a point processes follows from our assumption on $g_\varrho$). By $\overline{D}_m^{(\varrho)}$ we denote the $m$-th smallest point of $\xi_\varrho$ greater than zero and $\underline{D}_m^{(\varrho)}$ stands for the $m$-th largest point of $\xi_\varrho$ less than zero (with respect to the natural ordering). To neatly formulate a result about the asymptotic distributions of $\xi_\varrho$, $\overline{D}_m^{(\varrho)}$ and $\underline{D}_m^{(\varrho)}$, we use the following notation. For $\gamma>0$ and $a,b\in\RR$ with $a<b$ let us define 
$$\alpha_\varrho(a,b)=\frac{1}{2}\,\EE \sum_{(E,F)\in\eta_{t,\neq}^2} {\bf 1}\{m(E,F)\in K_\varrho,\,\varrho^{-\gamma} a < g_\varrho(E,F)\leq \varrho^{-\gamma} b\},$$
which is the expected number of pairs of flats with midpoint in $K_\varrho$ such that $\varrho^{-\gamma} a <g_\varrho(E,F)\leq\varrho^{-\gamma}b$. We further define
$$r_\varrho(u)=\sup_{E\in\AG^d_k}\Theta_t\big(\{F\in\AG^d_k: m(E,F)\in K_\varrho,\ -\varrho^{-\gamma}u \leq g_\varrho(E,F)\leq \varrho^{-\gamma}u\}\big)$$ for any $u>0$.
We are now in the position to formulate a two-sided version of Theorem 1.1 in \cite{SchulteThaele12}.

\begin{proposition}\label{prop:SchulteThaele12}
Let $\gamma>0$ and let $\nu$ be a $\sigma$-finite non-atomic Borel measure on $\RR$ such that
\begin{equation}\label{eq:conditionPointProcess}
\lim_{\varrho\rightarrow\infty}\alpha_\varrho(a,b)=\nu\big((a,b]\big)\quad\text{and}\quad\lim_{\varrho\rightarrow\infty}r_\varrho(u)=0
\end{equation}
for any $-\infty<a<b<\infty$ and $u>0$. Then there is a constant $C_u$ for every $u\geq 0$ such that
\begin{equation*}
\begin{split}
\Big|\PP\big(\varrho^\gamma \,\overline{D}_m^{(\varrho)}>u\big) & -e^{-\nu((0,u])}\sum_{i=0}^{m-1}\frac{\nu((0,u])^i}{i!}\Big| \\
& \leq \big|\nu((0,u])-\alpha_\varrho(0,u)\big|+C_u\sqrt{r_\varrho(u)}
\end{split}
\end{equation*}
and
\begin{equation*}
\begin{split}
\Big|\PP\big(\varrho^\gamma \underline{D}_m^{(\varrho)}< -u\big) & -e^{-\nu((-u,0])}\sum_{i=0}^{m-1}\frac{\nu((-u,0])^i}{i!}\Big|\\
& \leq \big|\nu((-u,0])-\alpha_\varrho(-u,0)\big|+C_u\sqrt{r_\varrho(u)}
\end{split}
\end{equation*}
for all $m\in\{1,2,3,\ldots\}$ and $\varrho\in[1,\infty)$. Furthermore, the rescaled point processes $\big(\varrho^\gamma \xi_\varrho\big)_{\varrho\geq 1}$ converge in distribution to a Poisson point process on $\RR$ with intensity measure $\nu$.
\end{proposition}

\begin{remark}\rm
In \cite{SchulteThaele12}, it is assumed that functions $(g_\varrho)_{\varrho\geq 1}$ are non-negative and that the measure $\nu$ satisfies $\nu(\dint u)=\beta\,\tau u^{\tau-1}\,{\bf 1}\{u>0\}\,\dint u$ for some constants $\beta,\tau>0$ (here, $\dint u$ stands for the element of the Lebesgue measure). However, these assumptions -- tailored to the applications in that paper -- can be relaxed so that Proposition \ref{prop:SchulteThaele12} can be shown by repeating literally the proof of Theorem 1.1 in \cite{SchulteThaele12}.\\
The assumptions on $\nu$ in the statement of Proposition \ref{prop:SchulteThaele12} ensure that a Poisson point process on $\RR$ with intensity measure $\nu$ exists; cf. Chapter 12 in \cite{Kallenberg}.\\
In contrast to $\alpha_\varrho(a,b)$, it is not necessary to consider $r_\varrho$ for arbitrary intervals $(a,b]$ because $(a,b]\subset [-u,u]$ for an appropriate choice of $u$ and $\lim\limits_{\varrho\rightarrow\infty}r_\varrho(u)=0$ for all $u>0$ already implies the same behaviour for all $(a,b]$ with $-\infty < a < b <\infty$.
\end{remark}

\subsection{Proof of Theorem \ref{thm:PointProcess}}
We apply Proposition \ref{prop:SchulteThaele12} to the functions $(g_\varrho)_{\varrho\geq 1}$ given by
$$g_\varrho(E,F)=\dist(E,F).$$ It remains to determine $\gamma$ and the measure $\nu$ as well as to check condition \eqref{eq:conditionPointProcess}. As a consequence of Theorem \ref{thm:Expectation}, we find that
\begin{eqnarray*}
\alpha_\varrho(a,b) & = & \EE\pi_t(K_\varrho,(\varrho^{-\gamma}b)_+)-\EE\pi_t(K_\varrho,(\varrho^{-\gamma}a)_+)\\
&=& {t^2\over 2}\kappa_{d-2k}V_d(K)\,\varrho^d\big((\varrho^{-\gamma}b)_+^{d-2k}-(\varrho^{-\gamma}a)_+^{d-2k}\big)\\
& & \qquad\qquad\qquad\qquad\qquad\times\int_{\LG_k^d}\int_{\LG_k^d}[M,L]\,\QQ(\dint L)\,\QQ(\dint M)
\end{eqnarray*}
for any reals $a<b$ (here $x_+=\max\{x,0\}$ for $x\in\RR$). Thus, choosing $\gamma=d/(d-2k)$ and putting $\nu$ as in the statement of the Theorem, we obtain 
$$\alpha_\varrho(a,b)=\nu\big((-\infty,b]\big)-\nu\big((-\infty,a]\big)=\nu\big((a,b]\big)$$
for all $-\infty<a<b<\infty$. Moreover, from \eqref{eq:proofCLT1} in Step 1 of the proof of Theorem \ref{thm:CLT} it follows that
\begin{eqnarray*}
r_\varrho(u) & = & \sup_{E\in\AG^d_k}\Theta_t\big(\{F\in\AG^d_k: m(E,F)\in K_\varrho,\,0\leq\dist(E,F)\leq\varrho^{-d/(d-2k)}u\}\big)\\
& \leq & t \kappa_k\kappa_{d-2k} (\diam(K)/2)^k u^{d-2k} \varrho^{-(d-k)},
\end{eqnarray*}
which because of $k<d$ tends to zero, as $\varrho\rightarrow\infty$ for any $u\geq 0$. So, Theorem \ref{thm:PointProcess} is a direct consequence of Proposition \ref{prop:SchulteThaele12}.\hfill $\Box$

\subsection{Proof of Theorem \ref{thm:Maximum}}
Let us apply Proposition \ref{prop:SchulteThaele12} to the family of functions $(g_\varrho)_{\varrho\geq 1}$ given by
$$g_\varrho(E,F)=\dist(E,F)-\sigma$$ so that the point process $\xi_\varrho^{(K,t)}$ in Theorem \ref{thm:Maximum} and the point process $\xi_\varrho$ in Proposition \ref{prop:SchulteThaele12} are related by $\xi_\varrho^{(K,t)}=\xi_\varrho+\sigma$.
As a consequence of Theorem \ref{thm:Expectation}, we find that in this case
\begin{eqnarray*}
\alpha_\varrho(a,b) & = & \EE\pi_t(K_\varrho,(\sigma+\varrho^{-\gamma}b)_+)-\EE\pi_t(K_\varrho,(\sigma+\varrho^{-\gamma}a)_+)\\
& = & \frac{t^2}{2} \kappa_{d-2k}V_d(K)\varrho^d\big((\sigma+\varrho^{-\gamma}b)_+^{d-2k}-(\sigma+\varrho^{-\gamma}a)_+^{d-2k}\,\big)\\
& & \qquad\qquad\qquad\qquad\qquad\times\int_{\LG^d_k}\int_{\LG^d_k} [M,L]\,\QQ(\dint L)\,\QQ(\dint M)
\end{eqnarray*}
for any reals $a<b$ (again, $x_+=\max\{x,0\}$ for $x\in\RR$). Since $\sigma+\varrho^{-\gamma}a\rightarrow\sigma$ and $\sigma+\varrho^{-\gamma}b\rightarrow\sigma$ as $\varrho\rightarrow\infty$ for all reals $a<b$, the measure $\nu$ is this time supported on the whole real axis. Together with $\gamma=d$ in the equation for $\alpha_\varrho(a,b)$ above, we obtain that $$\lim\limits_{\varrho\rightarrow\infty}\alpha_\varrho(a,b)=\beta\,(b-a)\quad {\rm for\;all\;reals}\quad a<b,$$ where
$\beta$ is given by \eqref{eq:betaMAX}. Moreover, there is a finite constant $C^{(1)}_{a,b}>0$ for any $a<b$ also depending on $K$, $t$ and $\sigma$ such that $\big|\alpha_\varrho(a,b)-\beta(b-a)\big|\leq C_{a,b}^{(1)} \varrho^{-d}$ for $\varrho\geq 1$. Using Lemma \ref{lem:Integral} and $[M,L]\leq 1$ in a similar way as in the proof of \eqref{eq:proofCLT1}, we have for any $u>0$,
\begin{equation*}
\begin{split}
& r_\varrho(u)\\
& = \sup_{E\in\AG^d_k}\Theta_t\big(\{F\in\AG^d_k: m(E,F)\in K_\varrho,\, -\varrho^{-d}u\leq \dist(E,F)-\sigma\leq \varrho^{-d}u\}\big)\\
& = t\sup_{M\in\LG^d_k,y\in M^\perp}\int_{\LG^d_k}\int_{(M+L)^\perp}V_k\big((K_\varrho-(x/2)-y)\cap M\big)\\
& \qquad\qquad\qquad\times\,{\bf 1}\{\,\sigma-\varrho^{-d}u\leq\|x\|\leq\sigma+\varrho^{-d}u\,\}\ell_{(M+L)^\perp}(\dint x) [M,L]\,\QQ(\dint L)\\
& \leq t\,\kappa_k(\diam(K)/2)^k \varrho^k \kappa_{d-2k}\left((\sigma+\varrho^{-d}u)^{d-2k}-(\sigma-\varrho^{-d}u)^{d-2k}\right)\\
& \leq C_u^{(2)} \varrho^{-(d-k)}
\end{split}
\end{equation*}
for all $\varrho\geq 1$ with a finite constant $C_u^{(2)}>0$ depending on $K$, $t$, $\sigma$ and $u$. Thus, the conditions in \eqref{eq:conditionPointProcess} are satisfied with $\gamma=d$ there and with $\nu$ equal to $\beta$ times the Lebesgue measure on $\RR$, where $\beta$ is given by \eqref{eq:betaMAX}. Whence Theorem \ref{thm:Maximum} is again a consequence of Proposition \ref{prop:SchulteThaele12}.\hfill $\Box$

\subsection{Proof of \eqref{eq:Maximum}}

Since each convex body includes a ball with positive radius, it is sufficient to assume that $K=B_r^d$, where $B_r^d\subset\RR^d$ is a ball with fixed radius $r>0$ around the origin. For $n=1,2,3,\ldots$ we define Poisson U-statistics
$$S_n=\frac{1}{2}\,\sum_{(E,F)\in\eta^2_{t,\neq}} {\bf 1}\{m(E,F)\in B_r^d,\, a_n<\dist(E,F)\leq b_n\}$$
with $a_n=2(3n-1)r$ and $b_n=6nr$.
From Theorem \ref{thm:Expectation}, it follows that
$$\EE S_n = \EE\pi_t(B_r^d,b_n)-\EE\pi_t(B_r^d,a_n) = c_1 \big(b_n^{d-2k}-a_n^{d-2k}\big)$$
with $c_1=\EE\pi_t(B_r^d,1)$. $S_n$ has a Wiener-It\^o chaos decomposition $$S_n=\EE S_n+I_1(f_1^{(n)})+I_2(f_2^{(n)})$$ with kernels
$$f_1^{(n)}(E)=\int_{\AG^d_k} {\bf 1}\{m(E,F)\in B_r^d,\, a_n<\dist(E,F)\leq b_n\}\, \Theta_t(\dint F)$$
and
$$f^{(n)}_2(E,F)=\frac{1}{2}\,{\bf 1}\{m(E,F)\in B_r^d,\, a_n<\dist(E,F)\leq b_n \}.$$
As a consequence of Lemma \ref{lem:Integral}, we have
$$f_1^{(n)}(E) \leq c_2 \left(b_n^{d-2k}-a_n^{d-2k}\right)\ \text{ for }\ E\in\AG^d_k,$$
where $c_2=t\kappa_k\kappa_{d-2k}\, r^k\int_{\LG^d_k}[M,L]\,\QQ(\dint L)$ and where $M\in\LG_k^d$ is $E$ shifted to the origin. Hence, we obtain
\begin{equation*}
\VV S_n = \|f_1^{(n)}\|_1^2 + 2 \|f_2^{(n)}\|_2^2 \leq 2c_2\,\big(b_n^{d-2k}-a_n^{d-2k}\big) \EE S_n + \EE S_n.
\end{equation*}
In order to belong to a pair $(E,F)$ with $a_n<\dist(E,F)\leq b_n$ and $m(E,F)\in B^d_r$, a flat $E\in \AG^d_k$ must satisfy $\frac{a_n}{2}-r<\dist(E,0)\leq \frac{b_n}{2}+r$. Since $\frac{a_n}{2}-r=(3n-2)r$ and $\frac{b_n}{2}+r=(3n+1)r$, the random variables $(S_n)_{n\in\NN}$ are determined by disjoint sets of $k$-flats and are independent by the Poisson assumption on $\eta_t$. As a consequence, the normalized random variables $\widetilde{S}_n=S_n/\EE S_n$ with $\EE\widetilde{S}_n=1$ for any $n\geq 1$ are independent, too. Together with the fact that $b_n^{d-2k}-a_n^{d-2k}\geq (b_n-a_n)^{d-2k}$, we obtain
\begin{equation*}
\VV\widetilde{S}_n = (\EE S_n)^{-2}\,\VV S_n \leq {2c_2\over c_1}+{1\over c_1\big(b_n^{d-2k}-a_n^{d-2k}\big)}\leq {2c_2\over c_1}+{1\over c_1(2r)^{d-2k}}<\infty.
\end{equation*}
Now, a version of the strong law of large numbers for independent, but not identically distributed random variables yields that
$$\lim_{N\rightarrow\infty} \frac{1}{N}\sum_{n=1}^N \widetilde{S}_n = 1\qquad{\rm with\ probability\ one;}$$ see \cite[Corollary 4.22]{Kallenberg}. Since each $S_n$ is almost surely bounded, this means that there is almost surely a sequence $(n_k)_{k\in\NN}$ with $S_{n_k}>0$ for all $k$. This implies \eqref{eq:Maximum}.
\hfill $\Box$

\subsection*{Acknowledgement}
We would like to thank Matthias Reitzner for valuable comments on a first version of the text.\\ This paper was written while the first author was visiting Case Western Reserve University. He was supported by the German Academic Exchange Service (DAAD).

\end{document}